\documentclass{article}
\usepackage{amssymb}
\usepackage{amsfonts}

\newtheorem{theorem}{Theorem}

\newtheorem{conjecture}[theorem]{Conjecture}
\newtheorem{corollary}[theorem]{Corollary}

\newtheorem{lemma}[theorem]{Lemma}

\newenvironment{proof}[1][Proof]{\noindent \textbf{#1.} }{\  \rule{0.5em}{0.5em}}
\input{tcilatex}
\begin{document}

\author{Zafer \c{S}iar \and Bing\"{o}l University, Department of
Mathematics, Bing\"{o}l/TURKEY \and zsiar@bingol.edu.tr}
\title{An exponential Diophantine equation related to the difference of
powers of two Fibonacci numbers}
\maketitle

\begin{abstract}
In this paper, we prove that there is no $x\geq 4$ such that the difference
of $x$-th powers of two consecutive Fibonacci numbers greater than $0$ is a
Lucas number. Also we show that the Diophantine equation 
\[
F_{n+l}^{x}-F_{n}^{x}=L_{r} 
\]%
with $l\in \left \{ 2,3,4\right \} ,$ $n>0,$ and $r\geq 0$ has no solutions
for $x\geq 4.$ Finally, we conjecture that the Diophantine equation 
\[
F_{n}^{x}-F_{m}^{x}=L_{r} 
\]%
with $(n,m)\neq (1,0),(2,0),$ and $r\geq 0$ has no solutions for $x\geq 4.$
\end{abstract}



\section{\protect\bigskip Introduction}

Let $(F_{n})$ and $(L_{n})$ be the sequences of Fibonacci numbers and of
Lucas numbers defined by $F_{0}=0,~F_{1}=1,$ $F_{n}=F_{n-1}+F_{n-2}$ and $%
L_{0}=2,L_{1}=1,~L_{n}=L_{n-1}+L_{n-2}$ for $n\geq 2,$ respectively. Binet
formulas for these numbers are 
\[
F_{n}=\frac{\alpha ^{n}-\beta ^{n}}{\sqrt{5}}\text{ and }L_{n}\text{ }%
=\alpha ^{n}+\beta ^{n}, 
\]%
where $\alpha =\dfrac{1+\sqrt{5}}{2}$ and $\beta =\dfrac{1-\sqrt{5}}{2},$
which are the roots of the characteristic equation $x^{2}-x-1=0.$ It can be
seen that $1<\alpha <2,$ $-1<\beta <0$ and $\alpha \beta =-1.$ The most
known identity related to these numbers is 
\begin{equation}
L_{n}=F_{n-1}+F_{n+1}.  \label{2.8}
\end{equation}%
If $n\geq 1$, then the relation between $F_{n}$ and $\alpha $ is given by 
\begin{equation}
\alpha ^{n-2}\leq F_{n}\leq \alpha ^{n-1}  \label{1.1}
\end{equation}%
and\ similarly, the relation between $n$-th Lucas number $L_{n}$ and $\alpha 
$ is \ 
\begin{equation}
\alpha ^{n-1}\leq L_{n}<2\alpha ^{n}.  \label{1.2}
\end{equation}%
The inequalities (\ref{1.1}) and (\ref{1.2}) can be proved by induction. For
more information about the Fibonacci and Lucas sequences with their
applications, one can see \cite{Ko}.

The problem of finding the perfect powers in the Fibonacci sequence was a
classical problem that attracted much attention over the last decades. One
can consult \cite{cohn} for Fibonacci numbers that are a square or twice a
square, and \cite{Bgud,Bgud2,togbe} for the similar studies. In all these
works, the authors have used elementary methods, congruences, modular
approach, and linear forms in logarithms. But, in recent years, many
mathematicians started to use particularly linear forms in logarithms of
algebraic numbers in order to solve some Diophantine equations including
Fibonacci, Lucas, Pell, and balancing numbers. For example, in \cite{togbe},
Marques and Togbe showed that if $s\geq 1$ is an integer such that $%
F_{m}^{s}+F_{m+1}^{s}$ is a Fibonacci number for all sufficiently large $m,$
then $s\in \left\{ 1,2\right\} .$ Then, Luca and Oyono, in \cite{oyono},
solved completely this problem. That is, they proved that the equation $%
F_{m}^{s}+F_{m+1}^{s}=F_{n}$ has no solutions $(m,n,s)$ with $m\geq 2$ and $%
s\geq 3.$ After that, in \cite{Ruiz}, the authors extended this problem to
the $k-$generalized Fibonacci numbers. In \cite{Rihane1}, Rihane et al.
tackled the Diophantine equation 
\[
P_{n}^{x}+P_{n+1}^{x}=P_{m} 
\]%
and gave all the solutions of this equation in nonnegative integers $m,n,x.$
Same authors, in \cite{Rihane2}, proved that the Diophantine equation 
\[
B_{n+1}^{x}-B_{n}^{x}=B_{m} 
\]%
has the solutions $(m,n,x)=(2n+2,n,2),(1,0,x),(0,n,0)$ in nonnegative
integers $m,n,$ and $x.$

On the other hand, the relation 
\begin{equation}
F_{n+2}-F_{n-2}=L_{n}  \label{1.3}
\end{equation}%
is well known. Motivated by this equality and the above mentioned studies,
we present a new problem. We will try to answer the question such that when
does the difference of $x$-th powers of\ any two Fibonacci numbers become a
Lucas number? Clearly, a trivial solution of this question for $x=1$ is seen
immediately from (\ref{1.3}). In this study, we show that the Diophantine
equation 
\[
F_{n+l}^{x}-F_{n}^{x}=L_{r} 
\]%
with $l\in \left\{ 1,2,3,4\right\} ,$ $n>0,$ and $r\geq 0$ has no solutions
for $x\geq 4.$ Finally, we conjecture that the Diophantine equation 
\begin{equation}
F_{n}^{x}-F_{m}^{x}=L_{r}  \label{1.4}
\end{equation}%
with $(n,m)\neq (1,0),(2,0),$ and $r\geq 0$ has no solutions for $x\geq 4.$
Here, we will prove this conjecture for $n\leq 2m+4.$ But, the proof of this
conjecture for $n>2m+4$ is really difficult.

Now let us give some inequalities, which will be useful for the proof of our
main theorem. We observe that the inequality 
\[
\dfrac{F_{n-1}}{F_{n}}\leq \frac{2}{3} 
\]%
holds for all $n\geq 3.$ This implies that 
\begin{equation}
\dfrac{F_{n}}{F_{m}}\geq \frac{3}{2}  \label{1.5}
\end{equation}%
for $m<n$ and $n\geq 3.$ Also, it follows that 
\[
\left( \dfrac{F_{m}}{F_{n}}\right) ^{x}\leq \left( \dfrac{F_{m}}{F_{n}}%
\right) ^{2}=\frac{1}{\left( F_{n}/F_{m}\right) ^{2}}\leq \frac{4}{9} 
\]%
for $x\geq 2$ and $n\geq 3.$ And thus, 
\begin{equation}
1-\left( \dfrac{F_{m}}{F_{n}}\right) ^{x}\geq \frac{1}{2}.  \label{1.6}
\end{equation}

Our main theorem is\ \ \ \ \ \ \ \ \ \ \ \ \ \ \ \ \ \ \ \ \ \ \ \ \ \ \ \ \
\ \ \ \ \ \ \ \ \ \ \ \ \ \ \ \ \ \ \ \ \ \ \ \ \ \ \ \ \ \ \ \ \ \ \ \ \ \
\ \ \ \ \ \ \ \ \ \ \ \ \ \ \ \ \ \ \ \ \ \ \ \ \ \ \ \ \ \ \ \ \ \ \ \ \ \
\ \ \ \ \ \ \ \ \ \ \ \ \ \ \ \ \ \ \ \ \ \ \ \ \ \ \ \ \ \ \ \ \ \ \ \ \ \
\ \ \ \ \ \ \ \ \ \ \ \ \ \ \ \ \ \ \ \ \ \ \ \ \ \ \ \ \ \ \ \ \ \ \ \ \ \
\ \ \ \ \ \ \ \ \ \ \ \ \ \ \ \ \ \ \ \ \ \ \ \ \ \ \ \ \ \ \ \ \ \ \ \ \ \
\ \ \ \ \ \ \ \ \ \ \ \ \ \ \ \ \ \ \ \ \ \ \ \ \ \ \ \ \ \ \ \ \ \ \ \ \ \
\ \ \ \ \ \ \ \ \ \ \ \ \ \ \ \ \ \ \ \ \ \ \ \ \ \ \ \ \ \ \ \ \ \ \ \ \ \
\ \ \ \ \ \ \ \ \ \ \ \ \ \ \ \ \ \ \ \ \ \ \ \ \ \ \ \ \ \ \ \ \ \ \ \ \ \
\ \ \ \ \ \ \ \ \ \ \ \ \ \ \ 

\begin{theorem}
\label{T4} Let $x$ be a positive integer and let $n,m,r$ be non-negative
integers such that $n\leq 2m+4$ if $m\neq 0.$ Then all solutions $\left(
n,m,r,x\right) $ of the Diophantine equation $F_{n}^{x}-F_{m}^{x}=L_{r}$ are
the elements of the sets%
\[
\left \{ 
\begin{array}{c}
(1,0,1,x),(2,0,1,x),\left( 3,0,3,2\right) ,\left( 3,0,0,1\right) , \\ 
\left( 4,0,2,1\right) ,\left( 3,1,1,1\right) ,\left( 3,2,1,1\right) ,\left(
3,1,2,2\right) ,\left( 3,2,2,2\right) , \\ 
\left( 3,1,4,3\right) ,\left( 3,2,4,3\right) ,\left( 4,1,0,1\right) ,\left(
4,2,0,1\right) ,\left( 5,4,0,1\right) , \\ 
\left( 4,3,1,1\right) ,\left( 5,2,3,1\right) ,\left( 6,5,2,1\right) ,\left(
6,1,4,1\right) ,\left( 5,3,2,1\right)%
\end{array}%
\right \} 
\]%
or $\left \{ \left( k+1,k-3,k-1,1\right) \right \} $ with $k\geq 4.$
\end{theorem}

\section{Auxiliary results}

In order to solve some Diophantine equations as in (\ref{1.4}), many
mathematicians have used Baker's theory of lower bounds for a nonzero linear
form in logarithms of algebraic numbers. Since such bounds are of crucial
importance in effectively solving of the Diophantine equation (\ref{1.4}),
we start with recalling some basic notions from algebraic number theory.

Let $\eta $ be an algebraic number of degree $d$ with minimal polynomial 
\[
a_{0}x^{d}+a_{1}x^{d-1}+\cdots +a_{d}=a_{0}\dprod\limits_{i=1}^{d}\left(
x-\eta ^{(i)}\right) \in \mathbb{Z}[x], 
\]%
where the $a_{i}$'s are relatively prime integers with $a_{0}>0$ and $\eta
^{(i)}$'s are conjugates of $\eta .$ Then 
\begin{equation}
h(\eta )=\frac{1}{d}\left( \log a_{0}+\dsum\limits_{i=1}^{d}\log \left( \max
\left\{ |\eta ^{(i)}|,1\right\} \right) \right)  \label{2.1}
\end{equation}%
is called the logarithmic height of $\eta .$ In particularly, if $\eta =a/b$
is a rational number with $\gcd (a,b)=1$ and $b>1,$ then $h(\eta )=\log
\left( \max \left\{ |a|,b\right\} \right) .$

The following properties of logarithmic height are found in many works
stated in the references:

\begin{equation}
h(\eta \pm \gamma )\leq h(\eta )+h(\gamma )+\log 2,  \label{2.2}
\end{equation}%
\begin{equation}
h(\eta \gamma ^{\pm 1})\leq h(\eta )+h(\gamma ),  \label{2.3}
\end{equation}%
\begin{equation}
h(\eta ^{m})=|m|h(\eta ).  \label{2.4}
\end{equation}%
The following theorem is deduced from Corollary 2.3 of Matveev \cite{Mtv}
and provides a large upper bound for the subscripts in the equation (\ref%
{1.4}) (also see Theorem 9.4 in \cite{Bgud}).

\begin{theorem}
\label{T1} Assume that $\gamma _{1},\gamma _{2},...,\gamma _{t}$ are
positive real algebraic numbers in a real algebraic number field $\mathbb{K}$
of degree $D$, $b_{1},b_{2},...,b_{t}$ are rational integers, and 
\[
\Lambda :=\gamma _{1}^{b_{1}}\cdots \gamma _{t}^{b_{t}}-1 
\]%
is not zero. Then 
\[
|\Lambda |>\exp \left( -1.4\cdot 30^{t+3}\cdot t^{4.5}\cdot D^{2}(1+\log
D)(1+\log B)A_{1}A_{2}\cdots A_{t}\right) , 
\]%
where 
\[
B\geq \max \left \{ |b_{1}|,...,|b_{t}|\right \} , 
\]%
and $A_{i}\geq \max \left \{ Dh(\gamma _{i}),|\log \gamma
_{i}|,0.16\right
\} $ for all $i=1,...,t.$
\end{theorem}

The following lemma, which will be used in the main theorem, gives a
sufficient condition for a rational number to be a convergent of a given
real number.

\begin{lemma}
\label{L7}(\cite{yann})Let $\gamma $ be a real number. Any non-zero rational
number $\frac{a}{b}$ with 
\[
\left \vert \gamma -\frac{a}{b}\right \vert <\frac{1}{2b^{2}} 
\]%
is a convergent of $\gamma .$
\end{lemma}

\begin{lemma}
\label{L3}(\cite{weger}, Weger's Lemma 2.2) Let $a,u\in 
\mathbb{R}
$ and $0<$ $a<1$. If $\left \vert u\right \vert <a,$ then 
\[
\left \vert \log (1+u)\right \vert <\frac{-\log (1-a)}{a}\cdot \left \vert
u\right \vert 
\]%
and 
\[
\left \vert u\right \vert <\frac{a}{1-e^{-a}}\cdot \left \vert
e^{u}-1\right
\vert . 
\]
\end{lemma}

The following lemma can be found in \cite{peter}.

\begin{lemma}
\label{L4}If $F_{n}|L_{m},$ then $n\leq 4.$
\end{lemma}

The following two theorems are given in \cite{cohn} and \cite{zfr},
respectively.

\begin{theorem}
\label{T2} If $L_{n}=x^{2},$ then $n=1,\,3$ and if $L_{n}=2x^{2},$ then $%
n=0,6.$
\end{theorem}

\begin{theorem}
\label{T2.1}If $L_{n}=L_{m}x^{2}$ with $m\geq 2,$ then $n=m.$
\end{theorem}

The following two theorems are proved in \cite{zfr2}.

\begin{theorem}
\label{T7}Let $1\leq m\leq n.$ Then all solutions of the equation $%
F_{n}+F_{m}=F_{r}$ are 
\[
(n,m,r)=(n,n-1,n+1),(1,1,3),(2,1,3),(2,2,3),(3,1,4). 
\]
\end{theorem}

\begin{theorem}
\label{T3}All solutions $(n,m,r,k)$ of the Diophantine equation $%
F_{k}=F_{n}+F_{m}+F_{r}$ with $1\leq r\leq m\leq n$ are the elements of the
sets 
\[
\left \{ \left( n,n-2,n-3,n+1\right) ,\left( n,n,n-1,n+2\right) \right \} 
\]%
and 
\[
\left \{ 
\begin{array}{c}
\left( 1,1,1,4\right) ,\left( 4,1,1,5\right) ,\left( 4,2,2,5\right) ,\left(
5,3,1,6\right) , \\ 
\left( 2,1,1,4\right) ,\left( 2,2,2,4\right) ,\left( 3,3,1,5\right)%
\end{array}%
\right \} . 
\]
\end{theorem}

The following lemma can be deduced from Theorem $2$ given in \cite{Bravo}.

\begin{lemma}
\label{L5}The Diophantine equation $L_{n}=2^{x}-1$ for some nonnegative
integers $n,x$ has only the solutions $(n,x)=(1,1),(2,2),(4,3).$
\end{lemma}

\section{The proof of Theorem\ \protect\ref{T4}\ \ \ \ \ \ \ \ \ \ \ \ \ \ \
\ \ \ \ \ \ \ \ \ \ \ \ \ \ \ \ \ \ \ \ \ \ \ \ \ \ \ \ \ \ \ \ \ \ \ \ \ \
\ \ \ \ \ \ \ \ \ \ \ \ \ \ \ \ \ \ \ \ \ \ \ \ \ \ \ \ \ \ \ \ \ \ \ \ \ \
\ \ \ \ \ \ \ \ \ \ \ \ \ \ \ \ \ \ \ \ \ \ \ \ \ \ \ \ \ \ \ \ \ \ }

\proof%
Let $x>0$ be an integer and let $n,m,r$ be non-negative integers such that $%
n\leq 2m+4$ if $m\neq 0.$ Assume that $F_{n}^{x}-F_{m}^{x}=L_{r}.$ It is
clear that $n\neq m$ since $L_{r}\neq 0$ for all integer $r.$ Then, $n>m.$
If $m=0,$ then we have $F_{n}^{x}=L_{r}.$ Since $F_{n}\mid L_{r},$ it
follows that $n\leq 4$ by Lemma \ref{L4}. It is obvious that $n\neq 0.$ If $%
n=1$ or $n=2,$ then $(n,m,r,x)=(1,0,1,x),(2,0,1,x).$ If $n=3$ or $n=4,$ then
we have $L_{r}=\square $ or $L_{r}=2\square ,$ or $L_{r}=L_{2}\square .$ In
these cases, we get the solutions $(n,m,r,x)=\left( 3,0,3,2\right) ,\left(
3,0,0,1\right) ,\left( 4,0,2,1\right) $ by Theorems \ref{T2} and \ref{T2.1}.
Now let $m\geq 1.$ If $m=1$ or $m=2,$ then we see that $n\geq 3.$ Assume
that $n=3.$ Hence, we have the equation $L_{r}=2^{x}-1.$ Then by Lemma \ref%
{L5}, we get 
\[
(n,m,r,x)\in \left \{ \left( 3,1,1,1\right) ,\left( 3,2,1,1\right) ,\left(
3,1,2,2\right) ,\left( 3,2,2,2\right) ,\left( 3,1,4,3\right) ,\left(
3,2,4,3\right) \right \} . 
\]%
From now on, assume that $m\geq 1,$ $n\geq 4,$ and $x\geq 1.$ Now let $x=1.$
Then we have the equation $F_{n}-F_{m}=L_{r},$ i.e., $%
F_{n}=F_{m}+F_{r-1}+F_{r+1}$ by (\ref{2.8}). By Theorems \ref{T7} and \ref%
{T3}, we obtain 
\[
(n,m,r,x)\in \left \{ 
\begin{array}{c}
\left( 3,1,1,1\right) ,\left( 3,2,1,1\right) ,\left( 4,3,1,1\right) ,\left(
4,1,0,1\right) ,\left( 4,2,0,1\right) , \\ 
\left( 5,4,0,1\right) ,\left( 5,2,3,1\right) ,\left( 6,5,2,1\right) ,\left(
6,1,4,1\right) ,\left( 5,3,2,1\right)%
\end{array}%
\right \} , 
\]%
and $(n,m,r,x)\in \left \{ \left( k+1,k-3,k-1,1\right) \right \} $ with $%
k\geq 4.$ Thus, we can suppose that $x\geq 2.$ Since $n\geq 4$ and $n>m,$ it
follows that $5\leq 3^{x}-2^{x}\leq F_{n}^{x}-F_{m}^{x}=L_{r},$ which
implies that $r\geq 4.$ On the other hand, using (\ref{1.1}), (\ref{1.2})
and (\ref{1.6}), we get 
\begin{equation}
\alpha ^{r-1}\leq L_{r}=F_{n}^{x}-F_{m}^{x}\leq F_{n}^{x}\leq \alpha
^{(n-1)x}  \label{a}
\end{equation}%
and also,%
\begin{equation}
2\alpha ^{r}>L_{r}=F_{n}^{x}-F_{m}^{x}=F_{n}^{x}\left( 1-\left( \frac{F_{m}}{%
F_{n}}\right) ^{x}\right) \geq \frac{F_{n}^{x}}{2}\geq \frac{\alpha ^{(n-2)x}%
}{2}.  \label{b}
\end{equation}%
If we make necessary calculations by using the inequalities (\ref{a}) and (%
\ref{b}), we get 
\begin{equation}
(m-3)x\leq (n-3)x\leq r<nx\leq \left( 2m+4\right) x,  \label{3.1}
\end{equation}%
where we used the facts that $n\geq 4$ and $m<n\leq 2m+4.$ Rearranging the
equation $F_{n}^{x}-F_{m}^{x}=L_{r}$ as $F_{n}^{x}-\alpha
^{r}=F_{m}^{x}+\beta ^{r}$ and taking absolute values of both sides of last
equality, we get%
\begin{equation}
\left \vert F_{n}^{x}-\alpha ^{r}\right \vert \leq F_{m}^{x}+\left \vert
\beta \right \vert ^{r}.  \label{3.2}
\end{equation}%
Dividing both sides of (\ref{3.2}) by $F_{n}^{x}$ yields to 
\begin{equation}
\left \vert 1-F_{n}^{-x}\alpha ^{r}\right \vert \leq \frac{1}{\left(
F_{n}/F_{m}\right) ^{x}}+\frac{\left \vert \beta \right \vert ^{r}}{F_{n}^{x}%
}\leq \frac{1}{\left( 1.5\right) ^{x}}+\frac{1}{\alpha ^{(n-2)x}}<\frac{2}{%
(1.5)^{x}},  \label{3.3}
\end{equation}%
where we used the inequality (\ref{1.5}) and the fact that $\alpha
^{(n-2)x}\geq \left( \frac{3}{2}\right) ^{x}$ for $m\geq 1$ and $n\geq 4.$
Put 
\[
\Lambda _{1}:=1-F_{n}^{-x}\alpha ^{r}. 
\]%
If $\Lambda _{1}=0,$ then we get $F_{n}^{x}=\alpha ^{r},$ which is
impossible since $\alpha ^{r}$ is irrational for all positive integers $r.$
So $\Lambda _{1}\neq 0.$ Now,\ let us apply Theorem \ref{T1} with $\gamma
_{1}:=\alpha ,~\gamma _{2}:=F_{n}$ and $b_{1}:=-x,~b_{2}:=r.$ Note that the
numbers $\gamma _{1}$ and $\gamma _{2}$ are positive real numbers and
elements of the field $\mathbb{K}=\mathbb{Q}(\sqrt{5}).$ It is obvious that
the degree of the field $\mathbb{K}$ is $2.$ So $D=2.$\ Moreover, since 
\[
h(\gamma _{1})=h(\alpha )=\dfrac{\log \alpha }{2}=\dfrac{0.4812...}{2} 
\]%
and 
\[
h(\gamma _{2})=h(F_{n})=n\log \alpha 
\]%
by (\ref{2.3}), we can take $A_{1}:=0.5$ and$~A_{2}:=2n\log \alpha .$ Also,
it is obvious that $r\geq x$ by (\ref{3.1}) since $n\geq 4.$ Therefore, we
can take $B=r.$ Thus, taking into account the inequality (\ref{3.3}) and
using Theorem \ref{T1}, we obtain%
\[
2(1.5)^{x}>\left \vert \Lambda _{1}\right \vert >\exp \left( -1.4\cdot
30^{5}\cdot 2^{4.5}\cdot 2^{2}(1+\log 2)(1+\log r)\left( 0.5\right) 2n\log
\alpha \right) . 
\]%
Taking logarithms in the above inequality, we get 
\[
x\log (1.5)-\log 2<2.51\cdot 10^{9}\cdot (1+\log r)\cdot n. 
\]%
Thus, it follows that 
\begin{equation}
x<6.2\cdot 10^{9}\cdot (1+\log nx)\cdot n  \label{z}
\end{equation}%
by (\ref{3.1}). Now we assume that $n\leq 270.$ Then, the inequality (\ref{z}%
) gives us that $x<6.43\cdot 10^{13}.$

Let 
\[
z_{1}:=r\log \alpha -x\log F_{n} 
\]%
and $u=e^{z_{1}}-1.$ By considering the inequality (\ref{3.3}), we have 
\[
\left \vert u\right \vert =\left \vert e^{z_{1}}-1\right \vert <\frac{2}{%
(1.5)^{x}}<0.9 
\]%
for $x\geq 2.$ Choosing $a=0.9$ in Lemma \ref{L3}, we get 
\[
|z_{1}|=\left \vert \log (u+1)\right \vert <\frac{-\log (1-0.9)}{(0.9)}\cdot 
\frac{2}{(1.5)^{x}}<\frac{5.12}{(1.5)^{x}}. 
\]%
Hence, it follows that 
\[
0<\left \vert r\log \alpha -x\log F_{n}\right \vert <\frac{5.12}{(1.5)^{x}}. 
\]%
Dividing both sides of this inequality by $x\log \alpha ,$ we get 
\begin{equation}
0<\left \vert \frac{\log F_{n}}{\log \alpha }-\frac{r}{x}\right \vert <\frac{%
11}{x\cdot (1.5)^{x}}.  \label{3.5}
\end{equation}%
Now assume that $x\geq 103.$ Then it can be seen that 
\[
\dfrac{(1.5)^{x}}{22}>6.23\cdot 10^{16}>6.43\cdot 10^{13}>x, 
\]%
and so we have 
\[
\left \vert \frac{\log F_{n}}{\log \alpha }-\frac{r}{x}\right \vert <\frac{11%
}{x\cdot (1.5)^{x}}<\dfrac{1}{2x^{2}}. 
\]%
Lemma \ref{L7} tells us that the rational number $\frac{r}{x}$ is a
convergent to $\gamma =\frac{\log F_{n}}{\log \alpha }.$ Then, let $%
[a_{0},a_{1},a_{2},...]$ be the continued fraction of $\gamma $ and let $%
p_{k}/q_{k}$ be its $k$-th convergent. Assume that $\frac{r}{x}=\dfrac{p_{t}%
}{q_{t}}$ for some $t.$ Then we have $1.7\cdot 10^{23}>q_{34}>6.43\cdot
10^{13}>x$ for every $n\in \left[ 4,270\right] .$ Thus $t\in
\{0,1,2,...,33\}.$ Furthermore, $a_{M}=\max \{a_{i}|i=0,1,2,...,33\}=1598.$
From the known properties of continued fraction, we get%
\begin{eqnarray*}
\left \vert \frac{\log F_{n}}{\log \alpha }-\frac{r}{x}\right \vert &=&\left
\vert \gamma -\dfrac{p_{t}}{q_{t}}\right \vert =\dfrac{1}{(\gamma
_{t+1}q_{t}+q_{t-1})q_{t}} \\
&=&\dfrac{1}{(\gamma _{t+1}+\frac{q_{t-1}}{q_{t}})q_{t}^{2}}>\dfrac{1}{%
(a_{t+1}+2)q_{t}^{2}}>\dfrac{1}{(a_{M}+2)q_{t}^{2}}\geq \dfrac{1}{1600x^{2}},
\end{eqnarray*}%
where we have used the facts that $a_{t}=\left \lfloor \gamma
_{t}\right
\rfloor $ and $q_{t-1}<q_{t}.$ Thus, from (\ref{3.5}), we obtain 
\[
\frac{11}{x\cdot (1.5)^{x}}>\dfrac{1}{1600x^{2}}, 
\]%
that is,%
\[
\frac{1}{12.46\cdot 10^{16}}>\frac{11}{(1.5)^{x}}>\dfrac{1}{1600x}\geq \frac{%
1}{10.288\cdot 10^{16}}, 
\]%
a contradiction. Therefore $x\leq 102.$ Taking into account the inequality (%
\ref{3.1}), a quick computation with Mathematica gives us that the equation $%
F_{n}^{x}-F_{m}^{x}=L_{r}$ has no solutions for $n\in \left[ 4,270\right] $
and $x\in \left[ 2,102\right] .$ Since this completes the analysis in the
case $n\in \left[ 4,270\right] ,$ from now on, we can assume that $n>270.$
This implies that $m\geq 134.$ Since $n\leq 2m+4$, from (\ref{z}), we can
write 
\begin{equation}
x<6.2\cdot 10^{9}\cdot (1+\log \left( 2m+4\right) x)\cdot \left( 2m+4\right)
.  \label{3.4}
\end{equation}%
Here, since $\left( m+2\right) x\geq 6$ implies that $(1+\log \left(
2m+4\right) x)\leq 2\log \left( m+2\right) x,$ we can rewrite the inequality
(\ref{3.4}) as 
\begin{equation}
x<24.8\cdot 10^{9}\cdot \left( m+2\right) \log \left( \left( m+2\right)
x\right)  \label{3.6}
\end{equation}%
If $x\leq m+2,$ then we have an inequality, which is better than inequality (%
\ref{3.6}). So, we are through. Contrast to this, if $x>m+2,$ then (\ref{3.6}%
) yields to us that $x<49.6\cdot 10^{9}\cdot \left( m+2\right) \log x,$
which can be rearranged as 
\begin{equation}
\frac{x}{\log x}<49.6\cdot 10^{9}\cdot \left( m+2\right) .  \label{3.7}
\end{equation}%
Using the fact that 
\[
\text{if }A\geq 3\text{ and }\frac{x}{\log x}<A,\text{ then }x<2A\log A, 
\]%
we obtain 
\[
x<99.2\cdot 10^{9}\cdot \left( m+2\right) \log \left( 49.6\cdot 10^{9}\cdot
\left( m+2\right) \right) , 
\]%
or%
\begin{equation}
x<595.2\cdot 10^{9}\left( m+2\right) \log \left( m+2\right) .  \label{3.8}
\end{equation}%
Now, put $y:=\dfrac{x}{\alpha ^{2m}}.$ Then, since $m\geq 134,$ from the
inequality (\ref{3.8}), we get 
\begin{equation}
y<\dfrac{595.2\cdot 10^{9}\left( m+2\right) \log \left( m+2\right) }{\alpha
^{2m}}<\dfrac{1}{\alpha ^{m}}.  \label{3.9}
\end{equation}%
Particularly, note that $y<\dfrac{1}{\alpha ^{m}}\leq \alpha
^{-134}<10^{-28}.$ On the other hand, it can be seen that 
\[
F_{n}^{x}=\dfrac{\alpha ^{nx}}{5^{x/2}}\left( 1-\dfrac{(-1)^{n}}{\alpha ^{2n}%
}\right) ^{x} 
\]%
and 
\[
F_{m}^{x}=\dfrac{\alpha ^{mx}}{5^{x/2}}\left( 1-\dfrac{(-1)^{m}}{\alpha ^{2m}%
}\right) ^{x}. 
\]%
Furthermore, we have 
\begin{eqnarray*}
0 &<&\left( 1-\dfrac{1}{\alpha ^{2n}}\right) ^{x}<1<\left( 1+\dfrac{1}{%
\alpha ^{2n}}\right) ^{x}=1+\frac{x}{\alpha ^{2n}}+\frac{x(x-1)}{\alpha ^{4n}%
}+... \\
&<&e^{y}<1+2y
\end{eqnarray*}%
because $y<10^{-28}.$ If we write $m$ instead of $n$ in the above
inequality, it holds. Thus, we see that 
\begin{equation}
\max \left \{ \left \vert F_{n}^{x}-\dfrac{\alpha ^{nx}}{5^{x/2}}\right
\vert ,\left \vert F_{m}^{x}-\dfrac{\alpha ^{mx}}{5^{x/2}}\right \vert
\right \} <\frac{2y\alpha ^{nx}}{5^{x/2}}.  \label{3.10}
\end{equation}%
Let us rearrange the equation $F_{n}^{x}-F_{m}^{x}=L_{r}$ as 
\[
\alpha ^{r}-\dfrac{\alpha ^{nx}}{5^{x/2}}+\dfrac{\alpha ^{mx}}{5^{x/2}}%
=F_{n}^{x}-\dfrac{\alpha ^{nx}}{5^{x/2}}-F_{m}^{x}+\dfrac{\alpha ^{mx}}{%
5^{x/2}}-\beta ^{r}. 
\]%
Taking absolute values of both sides of the above equality and using (\ref%
{3.10}), we get 
\[
\left \vert \alpha ^{r}-\dfrac{\alpha ^{nx}}{5^{x/2}}\left( 1-\alpha
^{(m-n)x}\right) \right \vert \leq \left \vert F_{n}^{x}-\dfrac{\alpha ^{nx}%
}{5^{x/2}}\right \vert +\left \vert F_{m}^{x}-\dfrac{\alpha ^{mx}}{5^{x/2}}%
\right \vert +\left \vert \beta \right \vert ^{r}<\frac{4y\alpha ^{nx}}{%
5^{x/2}}+\left \vert \beta \right \vert ^{r}. 
\]%
Dividing both sides of this inequality by $\dfrac{\alpha ^{nx}}{5^{x/2}},$
we obtain 
\begin{equation}
\left \vert \alpha ^{r-nx}\cdot 5^{x/2}+\alpha ^{(m-n)x}-1\right \vert <%
\dfrac{4}{\alpha ^{m}}+\frac{1}{\alpha ^{r}}\cdot \left( \frac{\sqrt{5}}{%
\alpha ^{n}}\right) ^{x},  \label{3.11}
\end{equation}%
where we used the fact that $\alpha \beta =-1$ and $y<\dfrac{1}{\alpha ^{m}}%
. $ Since 
\[
\left \vert \alpha ^{r-nx}\cdot 5^{x/2}-1\right \vert \leq \left \vert
\alpha ^{r-nx}\cdot 5^{x/2}+\alpha ^{(m-n)x}-1\right \vert +\alpha
^{(m-n)x}, 
\]%
from (\ref{3.11}), we get 
\begin{equation}
\left \vert \alpha ^{r-nx}\cdot 5^{x/2}-1\right \vert <\dfrac{4}{\alpha ^{m}}%
+\frac{1}{\alpha ^{r}}\cdot \left( \frac{\sqrt{5}}{\alpha ^{n}}\right) ^{x}+%
\frac{1}{\alpha ^{(n-m)x}}<\dfrac{4}{\alpha ^{m}}+\frac{1.05}{\alpha ^{x}}%
\text{,}  \label{3.12}
\end{equation}%
where we used the fact that $\left( \frac{\sqrt{5}}{\alpha ^{n}}\right)
^{x}\leq 0.05$ for $n>270$ and $r\geq x.$ Put 
\[
\Lambda _{2}:=1-\alpha ^{r-nx}\cdot 5^{x/2}. 
\]%
If $\Lambda _{2}=0,$ then we see that $\alpha ^{2(nx-r)}=5^{x},$ which is
possible only when $nx=r$ since $5^{x}\in 
\mathbb{Z}
.$ This is impossible since $r<nx$ by (\ref{3.1}). Therefore $\Lambda
_{2}\neq 0.$ Also, 
\begin{equation}
\left \vert \Lambda _{2}\right \vert <\dfrac{4}{\alpha ^{m}}+\frac{1.05}{%
\alpha ^{x}}<\frac{1}{2}  \label{c}
\end{equation}%
since $m\geq 134$ and $x\geq 2.$ Thus $\alpha ^{r-nx}\cdot 5^{x/2}\in \left[ 
\frac{1}{2},\frac{3}{2}\right] .$ Particularly, making necessary
calculations, it is seen that 
\begin{equation}
1.659x-0.86<nx-r<1.7x+1.46.  \label{3.13}
\end{equation}%
Let $k_{1}=\min \left \{ m,x\right \} .$ Then we see from (\ref{3.12}) that 
\begin{equation}
\left \vert \alpha ^{r-nx}\cdot 5^{x/2}-1\right \vert <\frac{5.05}{\alpha
^{k_{1}}}.  \label{3.14}
\end{equation}%
Now, let us apply Theorem \ref{T1} to the inequality (\ref{3.14}). Take $%
\gamma _{1}:=\sqrt{5},~\gamma _{2}:=\alpha ,$ and $b_{1}:=x,~b_{2}:=r-nx.$
Observe that the numbers $\gamma _{1}$ and $~\gamma _{2}$ are positive real
numbers and belong to the field $K=Q(\sqrt{5}).$ Therefore $D=2.$ Also,
since $h(\gamma _{1})=\log \sqrt{5}=0.804...,$ and $h(\gamma _{2})=\dfrac{%
\log \alpha }{2}=\dfrac{0.4812...}{2}$ by (\ref{2.1}), we can take $%
A_{1}:=1.61$ and $~A_{2}:=0.5.$ Besides, from (\ref{3.13}), it is clear that 
$nx-r<\left( 2.5\right) x$ for $x\geq 2.$ So, we can take $B=\left(
2.5\right) x\geq \max \left \{ |x|,|r-nx|\right \} .$ Thus, Theorem \ref{T1}
tells us that%
\[
\frac{5.05}{\alpha ^{k_{1}}}>\left \vert \Lambda _{2}\right \vert >\exp
\left( -C(1+\log 2)(1+\log \left( 2.5\right) x)\left( 1.61\right) \left(
0.5\right) \right) 
\]%
or 
\begin{equation}
k_{1}\log \alpha -\log 5.05<C(1+\log 2)(1+\log \left( 2.5\right) x)\left(
1.61\right) \left( 0.5\right) ,  \label{3.15}
\end{equation}%
where $C=1.4\cdot 30^{5}\cdot 2^{4.5}\cdot 2^{2}.$ If $k_{1}=x,$ then a
computer search with Mathematica gives us that $x<2.46\cdot 10^{11}.$ If $%
k_{1}=m,$ then, by\ using(\ref{3.8}), we get{\footnotesize 
\begin{equation}
m\log \alpha -\log 5.05<C(1+\log 2)(1+\log \left( 2.5\right) +\log
(595.2\cdot 10^{9}\left( m+2\right) \log \left( m+2\right) ))\left(
1.61\right) \left( 0.5\right) .  \label{3.16}
\end{equation}%
} With the help of a program in Mathematica\textit{, }the inequality (\ref%
{3.16}) gives us that $m<5.18\cdot 10^{11}.$ Substituting this value of $m$
into (\ref{3.8}), we obtain $x<8.32\cdot 10^{24}.$

Now, let 
\[
z_{2}:=x\log \sqrt{5}-(nx-r)\log \alpha 
\]%
and $u:=e^{z_{2}}-1.$ Then $\left \vert u\right \vert =\left \vert
e^{z_{2}}-1\right \vert =\left \vert \Lambda _{2}\right \vert <\frac{1}{2}$
by (\ref{c}). Thus, taking $a=1/2$ in Lemma \ref{L7} and making necessary
calculations, we get 
\[
|z_{2}|=\left \vert \log (1+u)\right \vert <\frac{-\log (1-\frac{1}{2})}{1/2}%
\left \vert u\right \vert <1.4\left( \dfrac{4}{\alpha ^{m}}+\frac{1.05}{%
\alpha ^{x}}\right) . 
\]%
That is, 
\[
0<\left \vert x\log \sqrt{5}-(nx-r)\log \alpha \right \vert <1.4\left( 
\dfrac{4}{\alpha ^{m}}+\frac{1.05}{\alpha ^{x}}\right) . 
\]%
Dividing both sides of the above inequality by $x\log \alpha ,$ we obtain 
\begin{equation}
\left \vert \frac{\log \sqrt{5}}{\log \alpha }-\frac{(nx-r)}{x}\right \vert <%
\frac{2.91}{x}\left( \dfrac{4}{\alpha ^{m}}+\frac{1.05}{\alpha ^{x}}\right) .
\label{3.17}
\end{equation}%
Since $m\geq 134,$ it follows that $\alpha ^{m}\geq \alpha
^{134}>10^{28}>1000x.$ Now we suppose that $x>100.$ Then it can be seen that 
$\alpha ^{x}>1000x.$ Hence, we can rewrite (\ref{3.17}) as 
\[
\left \vert \frac{\log \sqrt{5}}{\log \alpha }-\frac{(nx-r)}{x}\right \vert <%
\frac{2.91}{x}\left( \dfrac{4}{1000x}+\frac{1.05}{1000x}\right) <\frac{1}{%
66x^{2}}. 
\]%
This implies by Lemma \ref{L7} that the rational number $\frac{(nx-r)}{x}$
is a convergent to $\gamma =\frac{\log \sqrt{5}}{\log \alpha }.$ Now let $%
[a_{0},a_{1},a_{2},...]$ be the continued fraction of $\gamma $ and let $%
p_{k}/q_{k}$ be its $k$-th convergent. Assume that $\frac{(nx-r)}{x}%
=p_{t}/q_{t}$ for some $t.$ Then we have $2\cdot 10^{26}>q_{48}>8.32\cdot
10^{24}>x.$ Thus $t\in \{0,1,2,...,47\}.$ Furthermore, $a_{M}=\max
\{a_{i}|i=0,1,2,...,47\}=29.$ From the known properties of continued
fraction, we get%
\[
\frac{1}{66x^{2}}>\left \vert \frac{\log \sqrt{5}}{\log \alpha }-\frac{(nx-r)%
}{x}\right \vert =\left \vert \gamma -\dfrac{p_{t}}{q_{t}}\right \vert >%
\dfrac{1}{(a_{M}+2)q_{t}^{2}}\geq \dfrac{1}{31x^{2}}, 
\]%
a contradiction. So, $x\leq 100.$ Then $x<m.$ Hence, from (\ref{3.12}), we
get 
\begin{equation}
\left \vert \alpha ^{r-nx}\cdot 5^{x/2}-1\right \vert <\frac{5.05}{\alpha
^{x}}.  \label{3.18}
\end{equation}%
From (\ref{3.13}), we know that 
\[
1.659x-0.86<nx-r<1.7x+1.46. 
\]%
Put $t=nx-r.$ We found that the inequality (\ref{3.18}) is not satisfied\
for all $x\in \lbrack 2,100]$ and $t\in \left[ \left \lfloor
1.659x-0.86\right \rfloor ,\left \lceil 1.7x+1.46\right \rceil \right] .$
Thus the proof is completed.

Thus, we can give the following result.

\begin{corollary}
Let $x>0$ be an integer and let $n,r$ be nonnegative integers. Then all the
solutions $(n,r,x)$ of the Diophantine equation $F_{n+l}^{x}-F_{n}^{x}=L_{r}$
with $l\in \left \{ 1,2,3,4\right \} $ are given by 
\begin{eqnarray*}
(n,r,x) &=&(0,1,x),\left( 2,1,1\right) ,\left( 2,2,2\right) ,\left(
2,4,3\right) ,\left( 4,0,1\right) ,\left( 3,1,1\right) ,\left( 5,2,1\right) 
\text{ if }l=1, \\
(n,r,x) &=&(0,1,x),\left( 1,1,1\right) ,\left( 1,2,2\right) ,\left(
1,4,3\right) ,\left( 2,0,1\right) ,\left( 3,2,1\right) \text{ if }l=2, \\
(n,r,x) &=&\left( 0,3,2\right) ,\left( 0,0,1\right) ,\left( 1,0,1\right)
,\left( 2,3,1\right) \text{ if }l=3, \\
(n,r,x) &=&\left( 0,2,1\right) \text{ if }l=4.
\end{eqnarray*}
\end{corollary}

As one can see from the above result, the Diophantine equation 
\[
F_{n+l}^{x}-F_{n}^{x}=L_{r} 
\]%
with $l\in \left \{ 1,2,3,4\right \} ,$ $n>0,$ and $r\geq 0$ has no
solutions for $x\geq 4.$ If we pay attention, this equation has solutions
only for $n\leq 5.$ From the equations obtained by substituting these values
of $n$ (except for $n=0$) into the last equation, the equations, which have
a solution are given as follows:

\begin{corollary}
The Diophantine equation $3^{x}-2^{x}=L_{r}$ in nonnegative integers $r,x$
has only the solution $(r,x)=\left( 1,1\right) .$
\end{corollary}

\begin{corollary}
The Diophantine equation $5^{x}-3^{x}=L_{r}$ in nonnegative integers $r,x$
has only the solution $(r,x)=\left( 0,1\right) .$
\end{corollary}

\begin{corollary}
The Diophantine equation $8^{x}-5^{x}=L_{r}$ in nonnegative integers $r,x$
has only the solution $(r,x)=\left( 2,1\right) .$
\end{corollary}

\begin{corollary}
The Diophantine equation $3^{x}-1=L_{r}$ in nonnegative integers $r,x$ has
only the solution $(r,x)=\left( 0,1\right) .$
\end{corollary}

\begin{corollary}
The Diophantine equation $5^{x}-2^{x}=L_{r}$ in nonnegative integers $r,x$
has only the solution $(r,x)=\left( 2,1\right) .$
\end{corollary}

\begin{corollary}
The Diophantine equation $5^{x}-1=L_{r}$ in nonnegative integers $r,x$ has
only the solution $(r,x)=\left( 3,1\right) .$
\end{corollary}

\section{Concluding Remarks}

We were not able to solve Diophantine equation (\ref{1.4}) for $n>2m+4.$ But
we conjecture that the Diophantine equation (\ref{1.4}) has no solutions in
nonnegative integers $n,m,r,$ and $x$ when $n>2m+4.$ We think the following
conjecture is true and a computer search with Mathematica enables us to give
it.

\begin{conjecture}
The Diophantine equation (\ref{1.4}) with $(n,m)\neq (1,0),(2,0)$ has no
solutions for $x\geq 4.$
\end{conjecture}

\bigskip

\bigskip

\end{document}